\documentclass[12pt, fleqn]{article} 

\usepackage{amsmath,amssymb,amsthm}

\usepackage[numbers]{natbib}
\usepackage[margin=0.75in,centering]{geometry}
\usepackage{graphicx,float,pdflscape}
\usepackage{enumitem}
\usepackage[utf8]{inputenc}
\usepackage[section]{placeins}

\usepackage{makecell,longtable,booktabs,array}
\usepackage[table]{xcolor} 
\definecolor{lightgray}{rgb}{0.9, 0.9, 0.9}

\usepackage{microtype}
\usepackage[colorlinks=true, linkcolor=blue, citecolor=blue, urlcolor=blue]{hyperref}
\usepackage[noabbrev]{cleveref}

\newtheorem{thm}{Theorem}[section]

\newtheorem{lem}{Lemma}[section]

\theoremstyle{definition}
\newtheorem{defn}{Definition}[section]
\newtheorem{example}{Example}[section]
\theoremstyle{remark}

\title{Globally Simple Heffter Arrays $H(n;k)$ with $k \equiv 2 \pmod{4}$
       \thanks{This research was funded by TÜBİTAK, Grant/Award Number: 124F360}}
\author{Erik Pelttari \thanks{Department of Mathematics, Koç University. Email: epelttari22@ku.edu.tr}
        \and E. Şule Yazıcı \thanks{Department of Mathematics, Koç University. Email: eyazici@ku.edu.tr}
        }

\date{\today}

\begin{document}

\maketitle 

\begin{abstract}
Heffter arrays are combinatorial structures that can be used to construct orthogonal cyclic cycle decompositions of complete graphs. A square Heffter array $H(n;k)$ is a partially filled $n \times n$ array with $k$ filled cells in each row and column. Its entries are distinct and nonzero, form a half-set of the group $\mathbb{Z}_{2nk+1}$, and sum to zero in every row and column. If the natural orderings of the rows and columns of an array are simple, the array is called globally simple. The globally simple property guarantees the orthogonal cycle decompositions constructed from the Heffter array consist of simple cycles. While existence results are known for most congruence classes of the parameters $n$ and $k$ modulo 4 for square globally simple Heffter arrays $H(n;k)$, the existence for $k \equiv 2 \pmod{4}$ has remained largely open. In this paper, we present a construction for globally simple integer Heffter arrays $H(n;k)$ covering all admissible values of $n$ when $k \equiv 2 \pmod{4}$. This settles the existence question across all congruence classes.

\vspace{1em}
\noindent\textbf{Keywords and MSC codes:} Heffter arrays, globally simple, cycle decompositions. 05B30, 05C51
\end{abstract}

\section{Introduction}

\subsection{Overview and Definitions}

Introduced by Archdeacon \cite{archdeacon2015}, Heffter arrays are combinatorial structures that naturally extend Heffter's original difference problem \cite{heffter1896}. Heffter's first difference problem sought to partition the set $\{1,2,\dots,nk\}$ into triples $\{x,y,z\}$ such that $x+y+z \equiv 0 \pmod{2nk+1}$. A Heffter array generalizes this problem by allowing for arbitrary subset sizes and introducing two distinct partitions represented by the rows and columns of the array. By requiring each row and column in the array to sum to zero modulo $2nk+1$, Heffter arrays yield an explicit construction for orthogonal cycle decompositions of complete graphs under certain conditions \cite{archdeacon2015, cavenagh2020, costa2018}. Our work focuses on the relationship between Heffter arrays and orthogonal cycle decompositions. For an overview of the history, developments, and variants of Heffter arrays, see the survey by Dinitz and Pasotti \cite{dinitz2024}.

\begin{defn}\label{def:heffter}
A \textbf{Heffter array} $H(m,n;s,k)$ is a partially filled $m \times n$ array with entries from $\mathbb{Z}_{2nk+1}$, where $ms = nk$, such that:
\begin{enumerate}[label=(\roman*)]
    \item Each row contains exactly $s$ filled cells and each column contains exactly $k$ filled cells.
    \item For every $x \in \mathbb{Z}_{2nk+1} \setminus \{0\}$, exactly one of $x$ or $-x$ appears in the array (that is, the entries form a \textbf{half-set} of $\mathbb{Z}_{2nk+1}$).
    \item The sum of the elements in every row and every column is $0 \pmod{2nk+1}$.
\end{enumerate}
When the array is square ($m=n$ and thus $s=k$), it is denoted $H(n;k)$. We index the rows and columns by the set $\{0, 1, \dots, n-1\}$ throughout this paper.
\end{defn}

\begin{example}\label{ex:H106}
An example of an $H(10; 6)$ with entries from $\mathbb{Z}_{121}$. The construction of this array is detailed in \Cref{sec:Hn6_construction}.

\begin{center}
\renewcommand{\arraystretch}{1.25} % Slightly increases row height for readability
\begin{tabular}{|c|c|c|c|c|c|c|c|c|c|}
\hline
$9$ & $10$ & $-41$ & $19$ & & & $-39$ & $42$ & & \\ \hline
& $-24$ & $53$ & $-16$ & $-36$ & & & $-21$ & $44$ & \\ \hline
& & $1$ & $2$ & $-43$ & $17$ & & & $-23$ & $46$ \\ \hline
$48$ & & & $-26$ & $55$ & $-14$ & $-38$ & & & $-25$ \\ \hline
$-27$ & $50$ & & & $3$ & $4$ & $-45$ & $15$ & & \\ \hline
& $-29$ & $52$ & & & $-28$ & $57$ & $-12$ & $-40$ & \\ \hline
& & $-31$ & $54$ & & & $5$ & $6$ & $-47$ & $13$ \\ \hline
$-32$ & & & $-33$ & $56$ & & & $-30$ & $59$ & $-20$ \\ \hline
$-49$ & $11$ & & & $-35$ & $58$ & & & $7$ & $8$ \\ \hline
$51$ & $-18$ & $-34$ & & & $-37$ & $60$ & & & $-22$ \\ \hline
\end{tabular}
\end{center}
\end{example}

\begin{defn}\label{def:integer}
If the elements in every row and column sum to zero over the integers $\mathbb{Z}$ (rather than modulo $2nk+1$), the array is called an \textbf{integer Heffter array} \cite{archdeacon2_2015}. Note that \Cref{ex:H106} is an integer Heffter array.
\end{defn}

\begin{defn}\label{def:support}
In this context, we define the \textbf{support} of an array as the set of the absolute values of its filled cells. By condition (ii) of \Cref{def:heffter}, it follows that the support of a Heffter array is precisely the set $\{1, 2, \dots, nk\}$ \cite{archdeacon2015}. 
\end{defn}

In his foundational paper, Archdeacon utilized Heffter arrays to construct orthogonal cycle decompositions of complete graphs. We formalize this process with the following definitions.

\begin{defn}\label{def:orthogonal}
A cycle decomposition of a graph $G$ is a set of cycles, all of the same length, whose union contains every edge of $G$ exactly once. Two such decompositions, $\mathcal{C}_1$ and $\mathcal{C}_2$, are \textbf{orthogonal} if any cycle from $\mathcal{C}_1$ and any cycle from $\mathcal{C}_2$ intersect in at most one edge \cite{caro1999}.
\end{defn}

\begin{defn}\label{def:simple}
Let $T$ be a finite subset of an additive group $G$. Given an ordering $\omega = (t_1, t_2, \dots, t_k)$ of the elements in $T$, let $\displaystyle s_i = \sum_{j=1}^i t_j$ be the $i$-th partial sum of $\omega$. The ordering $\omega$ is said to be \textbf{simple} if its partial sums $s_1, s_2, \dots, s_k$ are pairwise distinct.   
\end{defn}

\begin{defn}\label{def:globsimp}
In a Heffter array $H(m,n;s,k)$, the filled cells of any row or column form a finite subset of the additive group $\mathbb{Z}_{2nk+1}$. If there exists a simple ordering for the elements of every row and every column, the array is said to be \textbf{simple}. If a Heffter array is simple under the natural left-to-right ordering of its rows and the top-to-bottom ordering of its columns, it is called \textbf{globally simple} \cite{costa2018}.
\end{defn}

If a Heffter array is globally simple, the sequence of partial sums in a row or column corresponds directly to the vertices of a starter cycle in $K_{2nk+1}$. By shifting these starter cycles through each of the $2nk+1$ vertices in the graph, the rows yield a cyclic decomposition of the complete graph $K_{2nk+1}$ into $s$-cycles, while the columns yield an orthogonal cyclic decomposition into $k$-cycles \cite{archdeacon2015}. If the rows or columns of a Heffter array do not have a simple ordering, the array yields orthogonal circuit decompositions instead. Given an arbitrary set of integers, finding a simple ordering is not a trivial problem. It connects to sequenceable groups introduced by Gordon in 1961 \cite{gordon1961} and continues with ongoing research into partial sums \cite{archdeacon2016, bode2005, costa2026}. For this reason, Heffter array constructions that are globally simple by design are particularly useful.

Before describing our construction, we state the relevant existence theorems for square, integer, and globally simple Heffter arrays. A square Heffter array $H(n;k)$ exists if and only if $3 \le k \le n$. For an integer Heffter array $H(n;k)$, existence further requires that $nk \equiv 0,3 \pmod{4}$ \cite{archdeacon2_2015, cavenagh2019, dinitz2017}. Finding Heffter arrays that are globally simple has been extensively studied across various congruence classes. For small values of $k$, Costa et al. constructed globally simple integer Heffter arrays for all valid $n, k$ pairs:

\begin{thm}\textup{\cite{costa2018}}
\label{thm:costa_small_k}
For $3 \le k \le 10$, a globally simple integer Heffter array $H(n;k)$ exists if and only if $k \le n$ and $nk \equiv 0,3 \pmod{4}$.
\end{thm}

For larger $k$, Burrage et al. successfully constructed globally simple integer Heffter arrays for the congruence classes where $k \equiv 0,3 \pmod{4}$.

\begin{thm}\textup{\cite{burrage2020}}
\label{thm:burrage_globally_simple}
For $3 \le k \le n$, a globally simple integer Heffter array $H(n;k)$ exists in the following cases:
\begin{enumerate}
    \item $k \equiv 0 \pmod{4}$,
    \item $n \equiv 1 \pmod{4}$ and $k \equiv 3 \pmod{4}$, and
    \item $n \equiv 0 \pmod{4}$ and $k \equiv 3 \pmod{4}$, provided $n \gg k$.
\end{enumerate}
\end{thm}

Recently, we established new results to resolve the $k \equiv 1 \pmod{4}$ class.

\begin{thm}\textup{\cite{pelttari2026}}
\label{thm:pelttari_globally_simple}
Let $k \le n$. For $k \equiv 1 \pmod{4}$ and $n \equiv 0,3 \pmod{4}$, there exists a globally simple integer Heffter array $H(n;k)$.
\end{thm}

These results establish the existence of globally simple Heffter arrays for three of the four congruence classes of $k$ modulo $4$, with the exception of sporadic open cases when $n \equiv 0 \pmod{4}$, $k \equiv 3 \pmod{4}$ and $n$ is close to $k$. Consequently, the primary remaining open problem is the existence of globally simple Heffter arrays for the final congruence class: $k \equiv 2 \pmod{4}$. For integer Heffter arrays, the condition $nk \equiv 0,3 \pmod{4}$ dictates that if $k \equiv 2 \pmod{4}$, then $n$ must be even. In this paper, we resolve this open problem by providing an explicit construction of integer Heffter arrays for $k \equiv 2 \pmod{4}$ with $k \ge 6$. This proves the following theorem.

\begin{thm}
\label{thm:nevenk2}
Let $6 \le k \le n$. For $k \equiv 2 \pmod{4}$ and $n$ even, there exists a globally simple integer Heffter array $H(n;k)$.
\end{thm}

In \Cref{sec:construction}, we define the explicit construction of the $H(n;k)$ arrays. In \Cref{sec:proof_simple}, we analyze the row and column partial sums to formally prove that the Heffter arrays we constructed are globally simple.

\section{Construction of Globally Simple $H(n;k)$ for $k \equiv 2 \pmod 4$} \label{sec:construction}

We construct a globally simple integer Heffter array $H(n;k)$ by initializing a base array $H(n;6)$ and augmenting it with sets of four diagonals until $k$ is reached. We index the rows and columns with the set $\{0,1,\dots,n-1\}$ throughout the construction. Operations on cell coordinates are performed modulo $n$.

\subsection{Diagonal Definition}\label{diag_def}

We define a diagonal of length $\ell$, using the notation from \cite{dinitz2017}. We denote by $\text{diag}(r, c, v, \delta, \delta_v, \ell)$ the set of $\ell$ filled cells starting at cell $(r, c)$ with initial value $v$. The coordinates and values of subsequent cells are determined by adding $\delta$ to the row and column indices and $\delta_v$ to the cell value. Formally, for $j \in \{0, 1, \dots, \ell-1\}$, the diagonal contains the entry $v + j\delta_v$ in cell $(r + j\delta, c + j\delta)$.

\subsection{Base Array $H(n;6)$ for Even $n$}\label{sec:Hn6_construction}

Although a globally simple integer $H(n;6)$ was constructed by Costa et al. \cite{costa2018}, our method requires a new $H(n;6)$ design to accommodate values of $k$ larger than $10$. The base array $H(n;6)$ is initialized using the elements positioned on the following ten diagonals. The support of these elements covers the integers $\{1, 2, \dots, 6n\}$. Note that \Cref{ex:H106_2} is constructed in this manner.

\begin{table}[h!]
\centering
\caption{Definitions and support of the diagonals for the base array $H(n;6)$.}
\label{tab:base_diagonals}
\footnotesize
\begin{tabular}{lll}
\toprule
\textbf{Diagonal} & \textbf{Definition} & \textbf{Support} \\
\midrule
$\mathbf{A}$ & $\text{diag}(2, 2, 1, 2, 2, n/2)$ & $\{1, 3, \dots, n-1\}$ \\
\addlinespace
$\mathbf{B}$ & $\text{diag}(n-1, n-1, -2n-2, 2, -2, n/2)$ & $\{2n+2, 2n+4, \dots, 3n\}$ \\
\addlinespace
$\mathbf{C}$ & $\text{diag}(2, 3, 2, 2, 2, n/2)$ & $\{2, 4, \dots, n\}$ \\
\addlinespace
$\mathbf{D}$ & $\text{diag}(n-1, 0, 5n+1, 2, 2, n/2)$ & $\{5n+1, 5n+3, \dots, 6n-1\}$ \\
\addlinespace
$\mathbf{E}$ & $\text{diag}(0, 2, -4n-1, 2, -2, n/2)$ & $\{4n+1, 4n+3, \dots, 5n-1\}$ \\
\addlinespace
$\mathbf{F}$ & $\text{diag}(n-3, n-1, -2n, 2, 2, n/2)$ & $\{2n, 2n-2, \dots, n+2\}$ \\
\addlinespace
$\mathbf{G}$ & $\text{diag}(0, 3, 2n-1, 2, -2, n/2)$ & $\{2n-1, 2n-3, \dots, n+1\}$ \\
\addlinespace
$\mathbf{H}$ & $\text{diag}(n-3, 0, -3n-2, 2, -2, n/2)$ & $\{3n+2, 3n+4, \dots, 4n\}$ \\
\addlinespace
$\mathbf{I}$ & $\text{diag}(0, n/2+2, 4n+2, 1, 2, n)$ & $\{4n+2, 4n+4, \dots, 6n\}$ \\
\addlinespace
$\mathbf{J}$ & $\text{diag}(1, n/2+2, -2n-1, 1, -2, n)$ & $\{2n+1, 2n+3, \dots, 4n-1\}$ \\
\bottomrule
\end{tabular}
\end{table}

\begin{example}\label{ex:H106_2}
An example of an $H(10; 6)$ with entries from $\mathbb{Z}_{121}$. The starting cell of each diagonal is marked with its respective subscript in the array. The cells covered by diagonals $\mathbf{A}$ and $\mathbf{I}$ are highlighted. 

\begin{center}
\renewcommand{\arraystretch}{1.25} 
% *{10} repeats the fixed-width centered column definition 10 times
\begin{tabular}{|*{10}{>{\centering\arraybackslash}p{1.1cm}|}}
\hline
\cellcolor{black!15}$9$ & $10$ & $-41_{\mathbf{E}}$ & $19_{\mathbf{G}}$ & & & $-39$ & \cellcolor{black!15}$42_{\mathbf{I}}$ & & \\ \hline
& $-24$ & $53$ & $-16$ & $-36$ & & & $-21_{\mathbf{J}}$ & \cellcolor{black!15}$44$ & \\ \hline
& & \cellcolor{black!15}$1_{\mathbf{A}}$ & $2_{\mathbf{C}}$ & $-43$ & $17$ & & & $-23$ & \cellcolor{black!15}$46$ \\ \hline
\cellcolor{black!15}$48$ & & & $-26$ & $55$ & $-14$ & $-38$ & & & $-25$ \\ \hline
$-27$ & \cellcolor{black!15}$50$ & & & \cellcolor{black!15}$3$ & $4$ & $-45$ & $15$ & & \\ \hline
& $-29$ & \cellcolor{black!15}$52$ & & & $-28$ & $57$ & $-12$ & $-40$ & \\ \hline
& & $-31$ & \cellcolor{black!15}$54$ & & & \cellcolor{black!15}$5$ & $6$ & $-47$ & $13$ \\ \hline
$-32_{\mathbf{H}}$ & & & $-33$ & \cellcolor{black!15}$56$ & & & $-30$ & $59$ & $-20_{\mathbf{F}}$ \\ \hline
$-49$ & $11$ & & & $-35$ & \cellcolor{black!15}$58$ & & & \cellcolor{black!15}$7$ & $8$ \\ \hline
$51_{\mathbf{D}}$ & $-18$ & $-34$ & & & $-37$ & \cellcolor{black!15}$60$ & & & $-22_{\mathbf{B}}$ \\ \hline
\end{tabular}
\end{center}
\end{example}

\subsection{Row and Column Sums}

To prove that this construction produces an integer Heffter array $H(n;k)$, we first show that every row and column in the base array $H(n;6)$ sums to zero.

\subsubsection{Row Sums of $H(n;6)$}

Each row $r \in \{0, 1, \dots, n-1\}$ contains exactly 6 filled cells. The diagonals intersecting the filled cells depend on the parity of $r$. \Cref{tab:row_sums} summarizes the entries of each diagonal across all row cases. Note that these diagonals are defined in \Cref{sec:Hn6_construction}.

\begin{table}[h!]
\centering
\caption{Diagonal entries and row sums for the base array $H(n;6)$.}
\label{tab:row_sums}
\small
\begin{tabular}{l c p{7.5cm} c}
\toprule
\textbf{Rows} & \textbf{Diagonals} & \textbf{Entries} & \textbf{Sum} \\
\midrule
%\multicolumn{4}{l}{\textit{Even Rows}} \\
$r \in \{2, 4, \dots, n-2\}$ & A, C, E, G, I, J & $(r-1)$, $r$, $(-4n-1-r)$, \newline $(2n-1-r)$, $(4n+2+2r)$, $(-2n+1-2r)$ & $0$ \\
\addlinespace
$r = 0$ & A, C, E, G, I, J & $(n-1)$, $n$, $(-4n-1)$, \newline $(2n-1)$, $(4n+2)$, $(-4n+1)$ & $0$ \\
\midrule
%\multicolumn{4}{l}{\textit{Odd Rows}} \\
$r \in \{1, 3, \dots, n-5\}$ & B, D, F, H, I, J & $(-2n-3-r)$, $(5n+2+r)$, $(-2n+3+r)$, \newline $(-3n-5-r)$, $(4n+2+2r)$, $(-2n+1-2r)$ & $0$ \\
\addlinespace
$r = n-3$ & B, D, F, H, I, J & $(-3n)$, $(6n-1)$, $(-2n)$, \newline $(-3n-2)$, $(6n-4)$, $(-4n+7)$ & $0$ \\
\addlinespace
$r = n-1$ & B, D, F, H, I, J & $(-2n-2)$, $(5n+1)$, $(-2n+2)$, \newline $(-3n-4)$, $6n$, $(-4n+3)$ & $0$ \\
\bottomrule
\end{tabular}
\end{table}

Hence, all row sums of $H(n;6)$ equal zero.

\subsubsection{Column Sums of $H(n;6)$}

Since diagonals I and J both begin at column $n/2+2$, and have the same step size, their contribution to column $c$ is constant:
\[
v_I(c) + v_J(c) =  2n + 1.
\]
Combining this constant sum with the four remaining diagonal entries in each column yields the column evaluation summarized in \Cref{tab:col_sums}.

\begin{table}[h!]
\centering
\caption{Diagonal entries and column sums for the base array $H(n;6)$.}
\label{tab:col_sums}
\small
\begin{tabular}{l c p{7.5cm} c}
\toprule
\textbf{Col Range} & \textbf{Diagonals} & \textbf{Entry Values (including $v_I + v_J$)} & \textbf{Sum} \\
\midrule
$c \in \{2, 4, \dots, n-2\}$ & A, D, E, H, I+J & $(c-1)$, $(5n+1+c)$, $(-4n+1-c)$, \newline $(-3n-2-c)$, $(2n+1)$ & $0$ \\
\addlinespace
$c = 0$ & A, D, E, H, I+J & $(n-1)$, $(5n+1)$, $(-5n+1)$, \newline $(-3n-2)$, $(2n+1)$ & $0$ \\
\midrule
$c \in \{3, 5, \dots, n-3\}$ & B, C, F, G, I+J & $(-2n-3-c)$, $(c-1)$, $(-2n+1+c)$, \newline $(2n+2-c)$, $(2n+1)$ & $0$ \\
\addlinespace
$c = 1$ & B, C, F, G, I+J & $(-2n-4)$, $n$, $(-2n+2)$, \newline $(n+1)$, $(2n+1)$ & $0$ \\
\addlinespace
$c = n-1$ & B, C, F, G, I+J & $(-2n-2)$, $(n-2)$, $(-2n)$, \newline $(n+3)$, $(2n+1)$ & $0$ \\
\bottomrule
\end{tabular}
\end{table}

Thus, every column sum in $H(n;6)$ is zero.

\subsection{Augmentation with $k-6$ New Diagonals}

To reach the required number of filled cells in each row and column, denoted by $k$, we augment the base array with $t = (k-6)/4$ sets of four diagonals. For each $i \in \{0, 1, \dots, t-1\}$, we define the $i$-th set of four added diagonals in \Cref{tab:added_diagonals}.

\begin{table}[h!]
\centering
\caption{Definitions and support of the $i$-th set of four added diagonals.}
\label{tab:added_diagonals}
\footnotesize
\begin{tabular}{lll}
\toprule
\textbf{Diagonal} & \textbf{Definition} & \textbf{Support} \\
\midrule
$\mathbf{D_{1+2i}}$ & $\text{diag}(1+2i, 0, (6n+4ni+1), 1, 4, n)$ & $\{6n+4ni+1, 6n+4ni+5, \dots, 10n+4ni-3\}$ \\
\addlinespace
$\mathbf{D_{2+2i}}$ & $\text{diag}(2+2i, 0, -(6n+4ni+3), 1, -4, n)$ & $\{6n+4ni+3, 6n+4ni+7, \dots, 10n+4ni-1\}$ \\
\addlinespace
$\mathbf{D_{n/2+2i}}$ & $\text{diag}(1+2i, n/2+1, -(6n+4ni+2), 1, -4, n)$ & $\{6n+4ni+2, 6n+4ni+6, \dots, 10n+4ni-2\}$ \\
\addlinespace
$\mathbf{D_{n/2+1+2i}}$ & $\text{diag}(2+2i, n/2+1, (6n+4ni+4), 1, 4, n)$ & $\{6n+4ni+4, 6n+4ni+8, \dots, 10n+4ni\}$ \\
\bottomrule
\end{tabular}
\end{table}

Each set of four added diagonals increases the number of filled cells in each row and column by $4$. As shown in \Cref{tab:added_diagonals}, the supports of $\mathbf{D_{1+2i}}$, $\mathbf{D_{2+2i}}$, $\mathbf{D_{n/2+2i}}$, and $\mathbf{D_{n/2+1+2i}}$ are contained in distinct congruence classes modulo $4$ and partition $\{6n+4ni+1, \dots, 10n+4ni\}$. Taking the union over all $0 \le i \le t-1$ (where $t=(k-6)/4$), the added diagonals provide the support $\{6n+1, 6n+2, \dots, 6n+4nt\} = \{6n+1, 6n+2, \dots, nk\}$. Combined with the support from the base array $H(n;6)$, the constructed $H(n;k)$ has the total support $\{1, 2, \dots, nk\}$. See \Cref{ex:H2018}.

\subsubsection{Sums for Augmented Array $H(n;k)$}

Recall from \Cref{diag_def} that for any diagonal $\text{diag}(r_0, c_0, v, \delta, \delta_v, n)$, the row and column indices increment by $\delta$ and starting from $v$, the successive cell values increment by $\delta_v$ at each step.

For any given row $r$, its intersection with the $i$-th set of added diagonals in \Cref{tab:added_diagonals} produces four entries:
\begin{itemize}
    \item \textbf{Diagonals $\mathbf{D_{1+2i}}$ and $\mathbf{D_{n/2+2i}}$:} Because both diagonals start in row $1+2i$ and have row step size $\delta = 1$, they intersect row $r$ after the same number of steps. Since their value step sizes $\delta_v$ are $+4$ and $-4$ respectively, the value changes from the starting cell cancel each other out at every row. Thus, their combined contribution to the sum of any row $r$ is constant and equals the sum of their initial values:
    \[
    (6n+4ni+1) + (-(6n+4ni+2)) = -1.
    \]
    
    \item \textbf{Diagonals $\mathbf{D_{2+2i}}$ and $\mathbf{D_{n/2+1+2i}}$:} Both diagonals start in row $2+2i$ with row step size $\delta = 1$ and opposite value step sizes $\delta_v$ of $-4$ and $+4$ respectively. By the same reasoning, their combined contribution to the sum of any row $r$ is constant:
    \[
    -(6n+4ni+3) + (6n+4ni+4) = 1.
    \]
\end{itemize}

Consequently, the contribution of each set of four added diagonals to any row $r$ is $(-1) + 1 = 0$.

The same analysis applies to the columns. For any given column $c$:
\begin{itemize}
    \item \textbf{Diagonals $\mathbf{D_{1+2i}}$ and $\mathbf{D_{2+2i}}$:} Both diagonals start in column $0$ with column step size $\delta = 1$ and opposite value step sizes $+4$ and $-4$, yielding a constant sum of:
    \[
    (6n+4ni+1) + (-(6n+4ni+3)) = -2.
    \]
    \item \textbf{Diagonals $\mathbf{D_{n/2+2i}}$ and $\mathbf{D_{n/2+1+2i}}$:} Both diagonals start in column $n/2+1$ with column step size $\delta = 1$ and opposite value step sizes $-4$ and $+4$, yielding a constant sum of:
    \[
    -(6n+4ni+2) + (6n+4ni+4) = +2.
    \]
\end{itemize}

The contribution of each set of four added diagonals to any column $c$ is therefore $(-2) + 2 = 0$. 

Together with the fact that the base array $H(n;6)$ has row and column sums of $0$ (\Cref{tab:row_sums,tab:col_sums}), we conclude that the sums of all rows and columns in the resulting $H(n;k)$ array are strictly zero.

%\subsection{Example: Globally Simple $H(20;18)$}

\begin{example}
\label{ex:H2018}
A globally simple $H(20;18)$ where the entries come from $\mathbb{Z}_{721}$. The base $H(20;6)$ array is highlighted in light gray, while the starting cells for the $12$ added diagonals are highlighted in dark gray.

\begin{table}[h]
    \centering
    \caption{Heffter Array $H(20;18)$}
    \label{tab:H2018}
    \setlength{\tabcolsep}{4pt}
    \scriptsize
        \begin{tabular}{|*{20}{c|}} 
            \hline
\cellcolor{black!15}19 & \cellcolor{black!15}20 & \cellcolor{black!15}-81 & \cellcolor{black!15}39 &  & 340 & -342 & 268 & -270 & 196 & -198 & \cellcolor{black!15}-79 & \cellcolor{black!15}82 &  & -339 & 341 & -267 & 269 & -195 & 197 \\
\hline
\cellcolor{black!40}121 & \cellcolor{black!15}-44 & \cellcolor{black!15}103 & \cellcolor{black!15}-36 & \cellcolor{black!15}-66 &  & 344 & -346 & 272 & -274 & 200 & \cellcolor{black!40}-122 & \cellcolor{black!15}-41 & \cellcolor{black!15}84 &  & -343 & 345 & -271 & 273 & -199 \\
\hline
\cellcolor{black!40}-123 & 125 & \cellcolor{black!15}1 & \cellcolor{black!15}2 & \cellcolor{black!15}-83 & \cellcolor{black!15}37 &  & 348 & -350 & 276 & -278 & \cellcolor{black!40}124 & -126 & \cellcolor{black!15}-43 & \cellcolor{black!15}86 &  & -347 & 349 & -275 & 277 \\
\hline
\cellcolor{black!40}201 & -127 & 129 & \cellcolor{black!15}-46 & \cellcolor{black!15}105 & \cellcolor{black!15}-34 & \cellcolor{black!15}-68 &  & 352 & -354 & 280 & \cellcolor{black!40}-202 & 128 & -130 & \cellcolor{black!15}-45 & \cellcolor{black!15}88 &  & -351 & 353 & -279 \\
\hline
\cellcolor{black!40}-203 & 205 & -131 & 133 & \cellcolor{black!15}3 & \cellcolor{black!15}4 & \cellcolor{black!15}-85 & \cellcolor{black!15}35 &  & 356 & -358 & \cellcolor{black!40}204 & -206 & 132 & -134 & \cellcolor{black!15}-47 & \cellcolor{black!15}90 &  & -355 & 357 \\
\hline
\cellcolor{black!40}281 & -207 & 209 & -135 & 137 & \cellcolor{black!15}-48 & \cellcolor{black!15}107 & \cellcolor{black!15}-32 & \cellcolor{black!15}-70 &  & 360 & \cellcolor{black!40}-282 & 208 & -210 & 136 & -138 & \cellcolor{black!15}-49 & \cellcolor{black!15}92 &  & -359 \\
\hline
\cellcolor{black!40}-283 & 285 & -211 & 213 & -139 & 141 & \cellcolor{black!15}5 & \cellcolor{black!15}6 & \cellcolor{black!15}-87 & \cellcolor{black!15}33 &  & \cellcolor{black!40}284 & -286 & 212 & -214 & 140 & -142 & \cellcolor{black!15}-51 & \cellcolor{black!15}94 &  \\
\hline
 & -287 & 289 & -215 & 217 & -143 & 145 & \cellcolor{black!15}-50 & \cellcolor{black!15}109 & \cellcolor{black!15}-30 & \cellcolor{black!15}-72 &  & 288 & -290 & 216 & -218 & 144 & -146 & \cellcolor{black!15}-53 & \cellcolor{black!15}96 \\
\hline
\cellcolor{black!15}98 &  & -291 & 293 & -219 & 221 & -147 & 149 & \cellcolor{black!15}7 & \cellcolor{black!15}8 & \cellcolor{black!15}-89 & \cellcolor{black!15}31 &  & 292 & -294 & 220 & -222 & 148 & -150 & \cellcolor{black!15}-55 \\
\hline
\cellcolor{black!15}-57 & \cellcolor{black!15}100 &  & -295 & 297 & -223 & 225 & -151 & 153 & \cellcolor{black!15}-52 & \cellcolor{black!15}111 & \cellcolor{black!15}-28 & \cellcolor{black!15}-74 &  & 296 & -298 & 224 & -226 & 152 & -154 \\
\hline
-158 & \cellcolor{black!15}-59 & \cellcolor{black!15}102 &  & -299 & 301 & -227 & 229 & -155 & 157 & \cellcolor{black!15}9 & \cellcolor{black!15}10 & \cellcolor{black!15}-91 & \cellcolor{black!15}29 &  & 300 & -302 & 228 & -230 & 156 \\
\hline
160 & -162 & \cellcolor{black!15}-61 & \cellcolor{black!15}104 &  & -303 & 305 & -231 & 233 & -159 & 161 & \cellcolor{black!15}-54 & \cellcolor{black!15}113 & \cellcolor{black!15}-26 & \cellcolor{black!15}-76 &  & 304 & -306 & 232 & -234 \\
\hline
-238 & 164 & -166 & \cellcolor{black!15}-63 & \cellcolor{black!15}106 &  & -307 & 309 & -235 & 237 & -163 & 165 & \cellcolor{black!15}11 & \cellcolor{black!15}12 & \cellcolor{black!15}-93 & \cellcolor{black!15}27 &  & 308 & -310 & 236 \\
\hline
240 & -242 & 168 & -170 & \cellcolor{black!15}-65 & \cellcolor{black!15}108 &  & -311 & 313 & -239 & 241 & -167 & 169 & \cellcolor{black!15}-56 & \cellcolor{black!15}115 & \cellcolor{black!15}-24 & \cellcolor{black!15}-78 &  & 312 & -314 \\
\hline
-318 & 244 & -246 & 172 & -174 & \cellcolor{black!15}-67 & \cellcolor{black!15}110 &  & -315 & 317 & -243 & 245 & -171 & 173 & \cellcolor{black!15}13 & \cellcolor{black!15}14 & \cellcolor{black!15}-95 & \cellcolor{black!15}25 &  & 316 \\
\hline
320 & -322 & 248 & -250 & 176 & -178 & \cellcolor{black!15}-69 & \cellcolor{black!15}112 &  & -319 & 321 & -247 & 249 & -175 & 177 & \cellcolor{black!15}-58 & \cellcolor{black!15}117 & \cellcolor{black!15}-22 & \cellcolor{black!15}-80 &  \\
\hline
 & 324 & -326 & 252 & -254 & 180 & -182 & \cellcolor{black!15}-71 & \cellcolor{black!15}114 &  & -323 & 325 & -251 & 253 & -179 & 181 & \cellcolor{black!15}15 & \cellcolor{black!15}16 & \cellcolor{black!15}-97 & \cellcolor{black!15}23 \\
\hline
\cellcolor{black!15}-62 &  & 328 & -330 & 256 & -258 & 184 & -186 & \cellcolor{black!15}-73 & \cellcolor{black!15}116 &  & -327 & 329 & -255 & 257 & -183 & 185 & \cellcolor{black!15}-60 & \cellcolor{black!15}119 & \cellcolor{black!15}-40 \\
\hline
\cellcolor{black!15}-99 & \cellcolor{black!15}21 &  & 332 & -334 & 260 & -262 & 188 & -190 & \cellcolor{black!15}-75 & \cellcolor{black!15}118 &  & -331 & 333 & -259 & 261 & -187 & 189 & \cellcolor{black!15}17 & \cellcolor{black!15}18 \\
\hline
\cellcolor{black!15}101 & \cellcolor{black!15}-38 & \cellcolor{black!15}-64 &  & 336 & -338 & 264 & -266 & 192 & -194 & \cellcolor{black!15}-77 & \cellcolor{black!15}120 &  & -335 & 337 & -263 & 265 & -191 & 193 & \cellcolor{black!15}-42 \\
\hline
        \end{tabular}
\end{table}
\end{example}

\section{\texorpdfstring{$H(n;k)$ is Globally Simple for $n$ Even and $k \equiv 2 \pmod{4}$}{H(n;k) is Globally Simple for n Even and k = 2 mod 4}} \label{sec:proof_simple}

We demonstrate that the array $H(n;k)$ constructed in \Cref{sec:construction} is globally simple. Given a starting cell, the proof analyzes the partial sums for the columns and rows, evaluating them modulo $2nk+1$. We denote the $p$-th partial sum of row $r$ and column $c$ by $S_R(r, p)$ and $S_C(c, p)$, respectively, with $0 \le r,c \le n-1$ and $0 \le p \le k-1$.

\subsection{Columns are simple}

In order to calculate the column partial sums, we need to choose a starting cell. For a fixed column $c$, the calculations are simplified by starting in cell $(c+n/2,c)$ and summing from top-to-bottom wrapping around when row $n-1$ is reached (recall that operations on cell coordinates are performed modulo $n$). Since the partial sums of any ordered set are distinct if and only if the partial sums of any cyclic shift of the ordered set are distinct, our choice of starting cell will not change whether the partial sums are distinct or not \cite{archdeacon2015, costa2018}. 

\subsubsection{Example: Column Partial Sums for the Globally Simple $H(20;18)$}

First we calculate the partial sums of the globally simple $H(20;18)$ found in \Cref{tab:H2018}. We denote the partial sums in each column using $S_C(c,p)$ so in column $0$ we have $S_C(0,0) = -158 \equiv 563 \pmod{721}, \ S_C(0,1) = 2, \ \dots, \  S_C(0,17) = 0$. 

\begin{table}[h]
    \centering
    \caption{Column Partial Sums Modulo $2nk+1=721$ for the globally simple $H(20;18)$ (sums involving base array entries are shaded in gray).}
    \label{tab:H2018_Col_Sums}
    \setlength{\tabcolsep}{4pt}
    \scriptsize
    % Column label + spacer + 18 data columns
    \begin{tabular}{lc*{18}{c}}
    \toprule
    Col 0 && -158 $\equiv$ 563 & 2 & 485 & 4 & 407 & 6 & \cellcolor{lightgray} 665 & \cellcolor{lightgray} 566 & \cellcolor{lightgray} 667 & \cellcolor{lightgray} 686 & 86 & 684 & 164 & 682 & 242 & 680 & \cellcolor{lightgray} 57 & \cellcolor{lightgray} 0 \\
    Col 1 && 559 & 2 & 481 & 4 & 403 & 6 & \cellcolor{lightgray} 27 & \cellcolor{lightgray} 710 & \cellcolor{lightgray} 9 & \cellcolor{lightgray} 686 & 90 & 684 & 168 & 682 & 246 & 680 & \cellcolor{lightgray} 59 & \cellcolor{lightgray} 0 \\
    Col 2 && 555 & 2 & 477 & 4 & 399 & 6 & \cellcolor{lightgray} 663 & \cellcolor{lightgray} 582 & \cellcolor{lightgray} 685 & \cellcolor{lightgray} 686 & 94 & 684 & 172 & 682 & 250 & 680 & \cellcolor{lightgray} 61 & \cellcolor{lightgray} 0 \\
    Col 3 && 551 & 2 & 473 & 4 & 395 & 6 & \cellcolor{lightgray} 45 & \cellcolor{lightgray} 9 & \cellcolor{lightgray} 11 & \cellcolor{lightgray} 686 & 98 & 684 & 176 & 682 & 254 & 680 & \cellcolor{lightgray} 63 & \cellcolor{lightgray} 0 \\
    Col 4 && 547 & 2 & 469 & 4 & 391 & 6 & \cellcolor{lightgray} 661 & \cellcolor{lightgray} 578 & \cellcolor{lightgray} 683 & \cellcolor{lightgray} 686 & 102 & 684 & 180 & 682 & 258 & 680 & \cellcolor{lightgray} 65 & \cellcolor{lightgray} 0 \\
    Col 5 && 543 & 2 & 465 & 4 & 387 & 6 & \cellcolor{lightgray} 43 & \cellcolor{lightgray} 9 & \cellcolor{lightgray} 13 & \cellcolor{lightgray} 686 & 106 & 684 & 184 & 682 & 262 & 680 & \cellcolor{lightgray} 67 & \cellcolor{lightgray} 0 \\
    Col 6 && 539 & 2 & 461 & 4 & 383 & 6 & \cellcolor{lightgray} 659 & \cellcolor{lightgray} 574 & \cellcolor{lightgray} 681 & \cellcolor{lightgray} 686 & 110 & 684 & 188 & 682 & 266 & 680 & \cellcolor{lightgray} 69 & \cellcolor{lightgray} 0 \\
    Col 7 && 535 & 2 & 457 & 4 & 379 & 6 & \cellcolor{lightgray} 41 & \cellcolor{lightgray} 9 & \cellcolor{lightgray} 15 & \cellcolor{lightgray} 686 & 114 & 684 & 192 & 682 & 270 & 680 & \cellcolor{lightgray} 71 & \cellcolor{lightgray} 0 \\
    Col 8 && 531 & 2 & 453 & 4 & 375 & 6 & \cellcolor{lightgray} 657 & \cellcolor{lightgray} 570 & \cellcolor{lightgray} 679 & \cellcolor{lightgray} 686 & 118 & 684 & 196 & 682 & 274 & 680 & \cellcolor{lightgray} 73 & \cellcolor{lightgray} 0 \\
    Col 9 && 527 & 2 & 449 & 4 & 371 & 6 & \cellcolor{lightgray} 39 & \cellcolor{lightgray} 9 & \cellcolor{lightgray} 17 & \cellcolor{lightgray} 686 & 122 & 684 & 200 & 682 & 278 & 680 & \cellcolor{lightgray} 75 & \cellcolor{lightgray} 0 \\
    Col 10 && 523 & 2 & 445 & 4 & 367 & 6 & \cellcolor{lightgray} 655 & \cellcolor{lightgray} 566 & \cellcolor{lightgray} 677 & \cellcolor{lightgray} 686 & 126 & 684 & 204 & 682 & 282 & 680 & \cellcolor{lightgray} 77 & \cellcolor{lightgray} 0 \\
    Col 11 && 599 & 2 & 521 & 4 & 443 & 6 & \cellcolor{lightgray} 37 & \cellcolor{lightgray} 9 & \cellcolor{lightgray} 19 & \cellcolor{lightgray} 686 & 130 & 684 & 208 & 682 & 286 & 680 & \cellcolor{lightgray} 79 & \cellcolor{lightgray} 0 \\
    Col 12 && 595 & 2 & 517 & 4 & 439 & 6 & \cellcolor{lightgray} 653 & \cellcolor{lightgray} 562 & \cellcolor{lightgray} 675 & \cellcolor{lightgray} 686 & 134 & 684 & 212 & 682 & 290 & 680 & \cellcolor{lightgray} 41 & \cellcolor{lightgray} 0 \\
    Col 13 && 591 & 2 & 513 & 4 & 435 & 6 & \cellcolor{lightgray} 35 & \cellcolor{lightgray} 9 & \cellcolor{lightgray} 21 & \cellcolor{lightgray} 686 & 138 & 684 & 216 & 682 & 294 & 680 & \cellcolor{lightgray} 43 & \cellcolor{lightgray} 0 \\
    Col 14 && 587 & 2 & 509 & 4 & 431 & 6 & \cellcolor{lightgray} 651 & \cellcolor{lightgray} 558 & \cellcolor{lightgray} 673 & \cellcolor{lightgray} 686 & 142 & 684 & 220 & 682 & 298 & 680 & \cellcolor{lightgray} 45 & \cellcolor{lightgray} 0 \\
    Col 15 && 583 & 2 & 505 & 4 & 427 & 6 & \cellcolor{lightgray} 33 & \cellcolor{lightgray} 9 & \cellcolor{lightgray} 23 & \cellcolor{lightgray} 686 & 146 & 684 & 224 & 682 & 302 & 680 & \cellcolor{lightgray} 47 & \cellcolor{lightgray} 0 \\
    Col 16 && 579 & 2 & 501 & 4 & 423 & 6 & \cellcolor{lightgray} 649 & \cellcolor{lightgray} 554 & \cellcolor{lightgray} 671 & \cellcolor{lightgray} 686 & 150 & 684 & 228 & 682 & 306 & 680 & \cellcolor{lightgray} 49 & \cellcolor{lightgray} 0 \\
    Col 17 && 575 & 2 & 497 & 4 & 419 & 6 & \cellcolor{lightgray} 31 & \cellcolor{lightgray} 9 & \cellcolor{lightgray} 25 & \cellcolor{lightgray} 686 & 154 & 684 & 232 & 682 & 310 & 680 & \cellcolor{lightgray} 51 & \cellcolor{lightgray} 0 \\
    Col 18 && 571 & 2 & 493 & 4 & 415 & 6 & \cellcolor{lightgray} 647 & \cellcolor{lightgray} 550 & \cellcolor{lightgray} 669 & \cellcolor{lightgray} 686 & 158 & 684 & 236 & 682 & 314 & 680 & \cellcolor{lightgray} 53 & \cellcolor{lightgray} 0 \\
    Col 19 && 567 & 2 & 489 & 4 & 411 & 6 & \cellcolor{lightgray} 29 & \cellcolor{lightgray} 710 & \cellcolor{lightgray} 7 & \cellcolor{lightgray} 686 & 162 & 684 & 240 & 682 & 318 & 680 & \cellcolor{lightgray} 55 & \cellcolor{lightgray} 0 \\
    \bottomrule
    \end{tabular}
\end{table}

\subsubsection{Column Partial Sums Modulo $2nk+1$}

To prove that this construction produces a globally simple array in the general case, we first state several lemmas used to derive the column partial sums for specific sets of columns. Recall that $t = (k-6)/4$ denotes the number of sets of four added diagonals.

\begin{lem}\label{lem:colbasevals}
Let $H(n;6)$ be the base array constructed in \Cref{sec:construction}. The non-empty entries of this array are exactly those detailed in \Cref{tab:Hn6colentries}, positioned top-to-bottom such that in each column $c$, four consecutive entries start at $(c-3, c)$ and two consecutive entries start at $(c+n/2-2, c)$. Furthermore, for any column $0 \le c \le n-1$, the sum of these four consecutive entries is $-2n-1$ and the sum of the two consecutive entries is always $2n+1$.
\end{lem}

\begin{table}[H]
\centering
\caption{Base $H(n;6)$ Column Entries}
\label{tab:Hn6colentries}
\setlength{\tabcolsep}{4pt} 
\footnotesize
\begin{tabular}{lccccccc}
\toprule
\textbf{Column} & \multicolumn{4}{c}{\textbf{Four Entries}} && \multicolumn{2}{c}{\textbf{Two Entries}} \\
\midrule
Column 0 & $-3n-2$, & $-5n+1$, & $5n+1$, & $n-1$ && $5n-2$, & $-3n+3$ \\
\addlinespace
Column 1 & $n+1$, & $-2n+2$, & $n$, & $-2n-4$ && $5n$, & $-3n+1$ \\
\midrule
Even $2 \le c \le n/2+1$ & $-3n-c-2$, & $-4n-c+1$, & $5n+c+1$, & $c-1$ && $5n+2c-2$, & $-3n-2c+3$ \\
\addlinespace
Odd $3 \le c \le n/2+1$ & $2n-c+2$, & $-2n+c+1$, & $c-1$, & $-2n-c-3$ && $5n+2c-2$, & $-3n-2c+3$ \\
\midrule
Even $n/2+2 \le c \le n-2$ & $-3n-c-2$, & $-4n-c+1$, & $5n+c+1$, & $c-1$ && $3n+2c-2$, & $-n-2c+3$ \\
\addlinespace
Odd $n/2+2 \le c \le n-3$ & $2n-c+2$, & $-2n+c+1$, & $c-1$, & $-2n-c-3$ && $3n+2c-2$, & $-n-2c+3$ \\
\midrule
Column $n-1$ & $n+3$, & $-2n$, & $n-2$, & $-2n-2$ && $5n-4$, & $-3n+5$ \\
\bottomrule
\end{tabular}
\end{table}

\begin{lem}\label{lem:coladdedvals}
In column $c$, the entries in diagonals $\mathbf{D_{1+2i}}$ and $\mathbf{D_{2+2i}}$ for $0 \le i \le t-1$ are given by:
\begin{equation*}
  4ni + 6n + 4c + 1 \quad \text{and} \quad -4ni - 6n - 4c - 3.
\end{equation*}
Similarly, for diagonals $\mathbf{D_{n/2+2i}}$ and $\mathbf{D_{n/2+1+2i}}$, the entries are given by:
\begin{equation*}
  -4ni - 6n - 4((c-n/2-1) \bmod{n}) - 2 \quad \text{and} \quad 4ni + 6n + 4((c-n/2-1) \bmod{n}) + 4.
\end{equation*}
\end{lem}

In preparation for calculating the complete list of column partial sums $S_C(c, p)$, we first determine the list of partial sums generated by the diagonals $\mathbf{D_{n/2+2i}}$ and $\mathbf{D_{n/2+1+2i}}$. Because the column partial sum calculations begin with this specific pair of diagonals starting in cell $(c+n/2, c)$, this initial sequence accounts for the first $2t$ partial sums.

\begin{lem}\label{lem:colsums_Dh}
For any column $0 \le c \le n-1$, summing the entries of the diagonals $\mathbf{D_{n/2+2i}}$ and $\mathbf{D_{n/2+1+2i}}$ top-to-bottom starting in cell $(c+n/2, c)$ with $0 \le i \le t-1$ yields the first $2t$ column partial sums $S_C(c, 0), \dots, S_C(c, 2t-1)$ modulo $2nk+1$. 
For columns $0 \le c \le n/2$, the sequence is:
\begin{equation*}
2nk-8n-4c+3, \ 2, \ 2nk-12n-4c+5, \ 4, \ \dots, \ 2nk-4n(t+1)-4c+2t+1, \ 2t.
\end{equation*}
For columns $n/2+1 \le c \le n-1$, the sequence is:
\begin{equation*}
2nk-4n-4c+3, \ 2, \ 2nk-8n-4c+5, \ 4, \ \dots, \ 2nk-4nt-4c+2t+1, \ 2t.
\end{equation*}
\end{lem}

We calculate these partial sums by adding the entries in column $c$ defined in \Cref{lem:coladdedvals}. The calculation splits into two cases based on the column $c$, generating the partial sums shown above.

Similarly, we determine the partial sums generated by the diagonals $\mathbf{D_{1+2i}}$ and $\mathbf{D_{2+2i}}$, again using \Cref{lem:coladdedvals}. As these sums occur after the $2t$ entries from $\mathbf{D_{n/2+2i}}$ and $\mathbf{D_{n/2+1+2i}}$ and the four consecutive base array entries, they will be shifted by the preceding partial sum of these $2t+4$ values. By \Cref{lem:colsums_Dh}, the cumulative sum after the first $2t$ diagonal entries is exactly $2t$. Because the subsequent four consecutive base array entries always sum to $-2n-1$, as seen in \Cref{lem:colbasevals}, this preceding partial sum evaluates to a constant $S_C(c,2t+3) = 2t - 2n - 1 \equiv 2nk-2n+2t \pmod{2nk+1}$ for every column.

\begin{lem}\label{lem:colsums_D12}
With the preceding partial sum $S_C(c,2t+3) = 2nk-2n+2t$, for any column $0 \le c \le n-1$, summing the entries of the diagonals $\mathbf{D_{1+2i}}$ and $\mathbf{D_{2+2i}}$ top-to-bottom starting in cell $(c+1, c)$ with $0 \le i \le t-1$ yields the next $2t$ partial sums $S_C(c, 2t+4), \dots, S_C(c, 4t+3)$:
\begin{align*}
& 4n+4c+2t, \ 2nk-2n+2t-2, \ 8n+4c+2t-2, \\
& 2nk-2n+2t-4, \ \dots, \ 4tn+4c+2, \ 2nk-2n.
\end{align*}
\end{lem}

To evaluate this sequence, we take the preceding cumulative sum $S_C(c,2t+3)=2nk-2n+2t$ and sequentially add the entries defined in \Cref{lem:coladdedvals}. Adding an entry of $\mathbf{D_{1+2i}}$ yields the partial sum $2nk-2n+2t+(4i+6)n+4c-2i+1 \equiv (4i+4)n+4c+2t-2i \pmod{2nk+1}$. Adding an entry of $\mathbf{D_{2+2i}}$ brings the cumulative sum to $2nk-2n+2t-2i-2$. Evaluating this sequence from $i=0$ to $t-1$ directly results in the $2t$ partial sums above. When $i=t-1$ the final partial sum in this sequence is $2nk-2n$ for all columns. 

\begin{lem}\label{lem:colsums_final}
Recall from \Cref{lem:colbasevals} that the sum of the two consecutive entries is $2n+1$. If the preceding partial sum is $S_C(c, k-3) = 2nk-2n$, then adding the two remaining base array entries from \Cref{tab:Hn6colentries} yields the final two partial sums $S_C(c, k-2)$ and $S_C(c, k-1) \equiv 0 \pmod{2nk+1}$. The penultimate sum $S_C(c, k-2)$ is determined by the column $c$:
\begin{description}
    \item[Columns $0 \le c \le n/2+1$:] $3n + 2c - 3$
    \item[Columns $n/2+2 \le c \le n-1$:] $n + 2c - 3$
\end{description}
In all cases, adding the final base array entry brings the cumulative sum to $S_C(c, k-1) \equiv 0 \pmod{2nk+1}$, completing the column sum.
\end{lem}

With these lemmas, we now construct the complete partial sum sequence $S_C(c,p)$. We evaluate the terms of $S_C(c,p)$ by sequentially adding: the $2t$ entries from $\mathbf{D_{n/2+2i}}$ and $\mathbf{D_{n/2+1+2i}}$ in \Cref{lem:colsums_Dh}, the four consecutive base $H(n;6)$ entries (\Cref{tab:Hn6colentries}), the $2t$ entries from $\mathbf{D_{1+2i}}$ and $\mathbf{D_{2+2i}}$ in \Cref{lem:colsums_D12}, and finally the remaining two consecutive base array entries which complete the column sum. 

\begin{thm}\label{thm:col_distinct}
For an $H(n;k)$ constructed as in \Cref{sec:construction}, the $k$ partial sums $S_C(c, p)$ in any column $c$ are distinct modulo $2nk+1$. 
\end{thm}

\begin{proof}
By evaluating and ordering the sums generated in Lemmas \ref{lem:colsums_Dh}--\ref{lem:colsums_final} modulo $2nk+1$, we obtain the following strict orderings $0 < x_1 < x_2 < \dots < x_{k-1}$ for each column. Since the column determines the values in the filled cells, we partition the columns into seven cases: column $0$, column $1$, even and odd columns for $2 \le c \le n/2+1$, even and odd columns for $n/2+2 \le c \le n-2$, and column $n-1$.

\textbf{Column 0:}
\begin{align*}
0 &< 2 < 4 < \dots < 2t < 3n-3 < 4n+2t < 8n+2t-2 < \dots < 4tn+2 \\
  &< 2nk-4n(t+1)+2t+1 < \dots < 2nk-8n+3 < 2nk-8n+2t \\
  &< 2nk-3n+2t-1 < 2nk-3n+2t+1 \\
  &< 2nk-2n < 2nk-2n+2 < \dots < 2nk-2n+2t.
\end{align*}

\textbf{Column 1:}
\begin{align*}
0 &< 2 < 4 < \dots < 2t < 2t+3 < n+2t+1 < 3n-1 < 4n+2t+4 < \dots < 4tn+6 \\
  &< 2nk-4n(t+1)+2t-3 < \dots < 2nk-8n-1 \\
  &< 2nk-2n < 2nk-2n+2 < \dots < 2nk-2n+2t < 2nk-n+2t+4.
\end{align*}

\textbf{Even Columns $2 \le c \le n/2+1$:}
\begin{align*}
0 &< 2 < 4 < \dots < 2t < 3n+2c-3 < 4n+4c+2t < \dots < 4tn+4c+2 \\
  &< 2nk-4n(t+1)-4c+2t+1 < \dots < 2nk-8n-4c+3 < 2nk-7n-2c+2t \\
  &< 2nk-3n-c+2t-1 \\
  &< 2nk-2n < 2nk-2n+2 < \dots < 2nk-2n+2t.
\end{align*}

\textit{Note: The partial sum $2nk-2n-c+2t+1$ also appears in this set of even columns. Recall $t \ge 1$. While $2nk-3n-c+2t-1 < 2nk-2n-c+2t+1$ is always true, the exact position of $2nk-2n-c+2t+1$ relative to the sequence $2nk-2n < 2nk-2n+2 < \dots < 2nk-2n+2t$ changes depending on $c$. However, because $c$ is even, $2nk-2n-c+2t+1$ is always odd. Since all terms in the sequence $2nk-2n+2i$ are even, this particular partial sum is guaranteed to be distinct no matter where it falls in this range.}

\textbf{Odd Columns $3 \le c \le n/2+1$:}
\begin{align*}
0 &< 2 < 4 < \dots < 2t < 2t+3 < c+2t+2 < 2n-c+2t+2 < 3n+2c-3 \\
  &< 4n+4c+2t < \dots < 4tn+4c+2 \\
  &< 2nk-4n(t+1)-4c+2t+1 < \dots < 2nk-8n-4c+3 \\
  &< 2nk-2n < 2nk-2n+2 < \dots < 2nk-2n+2t.
\end{align*}

\textbf{Even Columns $n/2+2 \le c \le n-2$:}
\begin{align*}
0 &< 2 < 4 < \dots < 2t < n+2c-3 < 4n+4c+2t < \dots < 4tn+4c+2 \\
  &< 2nk-4nt-4c+2t+1 < \dots < 2nk-4n-4c+3 < 2nk-7n-2c+2t \\
  &< 2nk-3n-c+2t-1 < 2nk-2n-c+2t+1 \\
  &< 2nk-2n < 2nk-2n+2 < \dots < 2nk-2n+2t.
\end{align*}

\textbf{Odd Columns $n/2+2 \le c \le n-3$:}
\begin{align*}
0 &< 2 < 4 < \dots < 2t < 2t+3 < c+2t+2 < 2n-c+2t+2 < n+2c-3 \\
  &< 4n+4c+2t < \dots < 4tn+4c+2 \\
  &< 2nk-4nt-4c+2t+1 < \dots < 2nk-4n-4c+3 \\
  &< 2nk-2n < 2nk-2n+2 < \dots < 2nk-2n+2t.
\end{align*}

\textbf{Column $n-1$:}
\begin{align*}
0 &< 2 < 4 < \dots < 2t < 2t+1 < n+2t+3 < 3n-5 < 8n+2t-4 < \dots < 4n(t+1)-2 \\
  &< 2nk-4n(t+1)+2t+5 < \dots < 2nk-8n+7 \\
  &< 2nk-2n < 2nk-2n+2 < \dots < 2nk-2n+2t < 2nk-n+2t+4.
\end{align*}

Because the inequalities in the proof of \Cref{thm:col_distinct} arrange the $k$ partial sums of each column into a strictly increasing sequence within the interval $[0, 2nk]$, no two partial sums in the same column can be equivalent. Consequently, the partial sums for every column are pairwise distinct. Since we also established in \Cref{lem:colsums_final} that the final cumulative sum of every column evaluates to $0 \pmod{2nk+1}$, these two facts together conclude the proof that the columns of the constructed array are simple.
\end{proof}

\subsection{Rows are simple}

Similar to the columns, we calculate the row partial sums by choosing a starting cell. Since the choice of starting cell does not affect the distinctness of the cyclic sequence of partial sums, we evaluate the sums starting from the leftmost column in the added diagonals. For a fixed row $r$, this means we begin the calculations by starting in cell $(r,r-(k-6)/2) = (r,r-2t)$ and summing from left-to-right, wrapping around when we reach column $n-1$ (recall that operations on cell coordinates are performed modulo $n$).

\subsubsection{Example: Row Partial Sums for the Globally Simple $H(20;18)$}

First we calculate the partial sums of the globally simple $H(20;18)$ found in \Cref{tab:H2018}. We denote the partial sums in each row using $S_R(r,p)$. In row $0$, we start at cell $(0, 14)$ and have $S_R(0,0) = -339 \equiv 382 \pmod{721}, \ S_R(0,1) = 2, \ \dots, \ S_R(0,17) = 0$.

\begin{table}[h]
    \centering
    \caption{Row Partial Sums Modulo $2nk+1=721$ for $H(20;18)$ (sums involving base array entries are shaded in gray).}
    \label{tab:H2018_Row_Sums}
    \setlength{\tabcolsep}{4pt}
    \scriptsize
    % Row label + spacer + 18 data columns
    \begin{tabular}{lc*{18}{c}}
    \toprule
    Row 0 && -339 $\equiv$ 382 & 2 & 456 & 4 & 530 & 6 & \cellcolor{lightgray} 25 & \cellcolor{lightgray} 45 & \cellcolor{lightgray} 685 & \cellcolor{lightgray} 3 & 343 & 1 & 269 & 720 & 195 & 718 & \cellcolor{lightgray} 639 & \cellcolor{lightgray} 0 \\
    Row 1 && 378 & 2 & 452 & 4 & 526 & 647 & \cellcolor{lightgray} 603 & \cellcolor{lightgray} 706 & \cellcolor{lightgray} 670 & \cellcolor{lightgray} 604 & 227 & 602 & 153 & 600 & 79 & 678 & \cellcolor{lightgray} 637 & \cellcolor{lightgray} 0 \\
    Row 2 && 374 & 2 & 448 & 4 & 602 & 6 & \cellcolor{lightgray} 7 & \cellcolor{lightgray} 9 & \cellcolor{lightgray} 647 & \cellcolor{lightgray} 684 & 311 & 682 & 237 & 680 & 83 & 678 & \cellcolor{lightgray} 635 & \cellcolor{lightgray} 0 \\
    Row 3 && 370 & 2 & 444 & 645 & 518 & 647 & \cellcolor{lightgray} 601 & \cellcolor{lightgray} 706 & \cellcolor{lightgray} 672 & \cellcolor{lightgray} 604 & 235 & 602 & 161 & 680 & 87 & 678 & \cellcolor{lightgray} 633 & \cellcolor{lightgray} 0 \\
    Row 4 && 366 & 2 & 520 & 4 & 594 & 6 & \cellcolor{lightgray} 9 & \cellcolor{lightgray} 13 & \cellcolor{lightgray} 649 & \cellcolor{lightgray} 684 & 319 & 682 & 165 & 680 & 91 & 678 & \cellcolor{lightgray} 631 & \cellcolor{lightgray} 0 \\
    Row 5 && 362 & 643 & 436 & 645 & 510 & 647 & \cellcolor{lightgray} 599 & \cellcolor{lightgray} 706 & \cellcolor{lightgray} 674 & \cellcolor{lightgray} 604 & 243 & 682 & 169 & 680 & 95 & 678 & \cellcolor{lightgray} 629 & \cellcolor{lightgray} 0 \\
    Row 6 && 438 & 2 & 512 & 4 & 586 & 6 & \cellcolor{lightgray} 11 & \cellcolor{lightgray} 17 & \cellcolor{lightgray} 651 & \cellcolor{lightgray} 684 & 247 & 682 & 173 & 680 & 99 & 678 & \cellcolor{lightgray} 627 & \cellcolor{lightgray} 0 \\
    Row 7 && 434 & 2 & 508 & 4 & 582 & 6 & \cellcolor{lightgray} 677 & \cellcolor{lightgray} 65 & \cellcolor{lightgray} 35 & \cellcolor{lightgray} 684 & 251 & 682 & 177 & 680 & 103 & 678 & \cellcolor{lightgray} 625 & \cellcolor{lightgray} 0 \\
    Row 8 && 430 & 2 & 504 & 4 & 578 & 6 & \cellcolor{lightgray} 13 & \cellcolor{lightgray} 21 & \cellcolor{lightgray} 653 & \cellcolor{lightgray} 684 & 255 & 682 & 181 & 680 & 107 & 678 & \cellcolor{lightgray} 623 & \cellcolor{lightgray} 0 \\
    Row 9 && 426 & 2 & 500 & 4 & 574 & 6 & \cellcolor{lightgray} 675 & \cellcolor{lightgray} 65 & \cellcolor{lightgray} 37 & \cellcolor{lightgray} 684 & 259 & 682 & 185 & 680 & 111 & 678 & \cellcolor{lightgray} 621 & \cellcolor{lightgray} 0 \\
    Row 10 && 422 & 2 & 496 & 4 & 570 & 6 & \cellcolor{lightgray} 15 & \cellcolor{lightgray} 25 & \cellcolor{lightgray} 655 & \cellcolor{lightgray} 684 & 263 & 682 & 189 & 680 & 115 & 678 & \cellcolor{lightgray} 619 & \cellcolor{lightgray} 0 \\
    Row 11 && 418 & 2 & 492 & 4 & 566 & 6 & \cellcolor{lightgray} 673 & \cellcolor{lightgray} 65 & \cellcolor{lightgray} 39 & \cellcolor{lightgray} 684 & 267 & 682 & 193 & 680 & 119 & 678 & \cellcolor{lightgray} 617 & \cellcolor{lightgray} 0 \\
    Row 12 && 414 & 2 & 488 & 4 & 562 & 6 & \cellcolor{lightgray} 17 & \cellcolor{lightgray} 29 & \cellcolor{lightgray} 657 & \cellcolor{lightgray} 684 & 271 & 682 & 197 & 680 & 123 & 678 & \cellcolor{lightgray} 615 & \cellcolor{lightgray} 0 \\
    Row 13 && 410 & 2 & 484 & 4 & 558 & 6 & \cellcolor{lightgray} 671 & \cellcolor{lightgray} 65 & \cellcolor{lightgray} 41 & \cellcolor{lightgray} 684 & 275 & 682 & 201 & 680 & 127 & 678 & \cellcolor{lightgray} 613 & \cellcolor{lightgray} 0 \\
    Row 14 && 406 & 2 & 480 & 4 & 554 & 6 & \cellcolor{lightgray} 19 & \cellcolor{lightgray} 33 & \cellcolor{lightgray} 659 & \cellcolor{lightgray} 684 & 279 & 682 & 205 & 680 & 131 & 678 & \cellcolor{lightgray} 611 & \cellcolor{lightgray} 0 \\
    Row 15 && 402 & 2 & 476 & 4 & 550 & 6 & \cellcolor{lightgray} 669 & \cellcolor{lightgray} 65 & \cellcolor{lightgray} 43 & \cellcolor{lightgray} 684 & 283 & 682 & 209 & 680 & 135 & 678 & \cellcolor{lightgray} 609 & \cellcolor{lightgray} 0 \\
    Row 16 && 398 & 2 & 472 & 4 & 546 & 6 & \cellcolor{lightgray} 21 & \cellcolor{lightgray} 37 & \cellcolor{lightgray} 661 & \cellcolor{lightgray} 684 & 287 & 682 & 213 & 680 & 139 & 678 & \cellcolor{lightgray} 607 & \cellcolor{lightgray} 0 \\
    Row 17 && 394 & 2 & 468 & 4 & 542 & 6 & \cellcolor{lightgray} 667 & \cellcolor{lightgray} 65 & \cellcolor{lightgray} 25 & \cellcolor{lightgray} 684 & 291 & 682 & 217 & 680 & 143 & 678 & \cellcolor{lightgray} 605 & \cellcolor{lightgray} 0 \\
    Row 18 && 390 & 2 & 464 & 4 & 538 & 6 & \cellcolor{lightgray} 23 & \cellcolor{lightgray} 41 & \cellcolor{lightgray} 663 & \cellcolor{lightgray} 684 & 295 & 682 & 221 & 680 & 147 & 678 & \cellcolor{lightgray} 603 & \cellcolor{lightgray} 0 \\
    Row 19 && 386 & 2 & 460 & 4 & 534 & 6 & \cellcolor{lightgray} 685 & \cellcolor{lightgray} 65 & \cellcolor{lightgray} 27 & \cellcolor{lightgray} 684 & 299 & 682 & 225 & 680 & 151 & 678 & \cellcolor{lightgray} 601 & \cellcolor{lightgray} 0 \\
    \bottomrule
    \end{tabular}
\end{table}

\subsubsection{Row Partial Sums Modulo $2nk+1$}
As with the columns, we state several lemmas used to derive the row partial sums for specific sets of rows.

\begin{lem}\label{lem:rowbasevals}
Let $H(n;6)$ be the base array constructed in \Cref{sec:construction}. The non-empty entries of this array are exactly those detailed in \Cref{tab:Hn6rowentries}, positioned such that in each row $r$, four consecutive entries start at cell $(r, r)$ and two consecutive entries start at cell $(r, r+n/2+1)$. Furthermore, we see that the sum of the four consecutive entries is $-3$ for row $0$ and $-2n-3$ for all rows $1 \le r \le n-1$, while the sum of the two consecutive entries is $3$ for row $0$ and $2n+3$ for all rows $1 \le r \le n-1$.
\end{lem}

\begin{table}[h]
\centering
\caption{Base $H(n;6)$ Row Entries}
\label{tab:Hn6rowentries}
\setlength{\tabcolsep}{4pt} 
\footnotesize
\begin{tabular}{lccccccc}
\toprule
\textbf{Row} & \multicolumn{4}{c}{\textbf{Four Consecutive Entries}} && \multicolumn{2}{c}{\textbf{Two Consecutive Entries}} \\
\midrule
Row 0 & $n-1$, & $n$, & $-4n-1$, & $2n-1$ && $-4n+1$, & $4n+2$ \\
\addlinespace
Odd $1 \le r \le n-5$ & $-2n-r-3$, & $5n+r+2$, & $-2n+r+3$, & $-3n-r-5$ && $-2n-2r+1$, & $4n+2r+2$ \\
\addlinespace
Even $2 \le r \le n-2$ & $r-1$, & $r$, & $-4n-r-1$, & $2n-r-1$ && $-2n-2r+1$, & $4n+2r+2$ \\
\addlinespace
Row $n-3$ & $-3n$, & $6n-1$, & $-2n$, & $-3n-2$ && $-4n+7$, & $6n-4$ \\
\addlinespace
Row $n-1$ & $-2n-2$, & $5n+1$, & $-2n+2$, & $-3n-4$ && $-4n+3$, & $6n$ \\
\bottomrule
\end{tabular}
\end{table}

\begin{lem}\label{lem:rowaddedvals}
The entries in diagonals $\mathbf{D_{2+2i}}$ and $\mathbf{D_{1+2i}}$ for a fixed row $r$ with $0 \le i \le t-1$ are given by:

\begin{equation*}
  -6n-4ni-3-4((r-2i-2) \bmod{n}) \quad \text{and} \quad 6n+4ni+1+4((r-2i-1) \bmod{n}).
\end{equation*}

Similarly, the entries in row $r$ for diagonals $\mathbf{D_{n/2+1+2i}}$ and $\mathbf{D_{n/2+2i}}$ are given by:

\begin{equation*}
  6n+4ni+4+4((r-2i-2) \bmod{n}) \quad \text{and} \quad -6n-4ni-2-4((r-2i-1) \bmod{n}).
\end{equation*}

\end{lem}

\begin{lem}\label{lem:rowaddedsums_D12}
For any row $0 \le r \le n-1$, summing the entries in the diagonals $\mathbf{D_{2+2i}}$ and $\mathbf{D_{1+2i}}$ left-to-right starting in cell $(r, r-2t)$ with $0 \le i \le t-1$ yields the first $2t$ row partial sums $S_R(r, 0), \dots, S_R(r, 2t-1)$. The pair of partial sums modulo $2nk+1$ generated at step $i$ is determined by four cases based on the row $r$:
\begin{description}
    \item[Case 1: Row $0$.]
    The pair of sums is: $2nk-6n-4nt+8t-2+(4n-6)i \quad \text{and} \quad 2i+2$.
    
    \item[Case 2: Odd Rows $1 \le r \le 2t-1$.]
    For $i < t-(r+1)/2$, the pair is $2nk-6n-4nt+8t-4r-2+(4n-6)i$ and $2i+2$. For $i \ge t-(r+1)/2$, the pair is $2nk-6n-4nt+8t-4r-2+(4n-6)i$ and $2nk-4n+2i+3$.
    
    \item[Case 3: Even Rows $2 \le r \le 2t-2$.]
    For $i < t-r/2$, the pair is $2nk-6n-4nt+8t-4r-2+(4n-6)i$ and $2i+2$. For $i \ge t-r/2$, the pair is $2nk-2n-4nt+8t-4r-2+(4n-6)i$ and $2i+2$.
    
    \item[Case 4: Rows $2t \le r \le n-1$.]
    The pair of sums is: $2nk-2n-4nt+8t-4r-2+(4n-6)i \quad \text{and} \quad 2i+2$.
\end{description}
\end{lem}

Next, we determine the partial sums generated by the diagonals $\mathbf{D_{n/2+1+2i}}$ and $\mathbf{D_{n/2+2i}}$. As these sums occur after the first $2t$ added diagonal entries and the four consecutive base array entries, they are shifted by the preceding partial sum $S_R(r, 2t+3)$. By evaluating \Cref{lem:rowaddedsums_D12} at $i=t-1$ and adding the sum of the four consecutive base entries from \Cref{lem:rowbasevals}, we obtain three values for $S_R(r, 2t+3)$ modulo $2nk+1$:
\begin{itemize}
    \item \textbf{Row $0$:} $S_R(0, 2t+3) = 2t - 3$ (or $2nk$ when $t=1$).
    \item \textbf{Odd Rows $1 \le r \le 2t-1$:} $S_R(r, 2t+3) = 2t - 6n - 3 \equiv 2nk - 6n + 2t - 2 \pmod{2nk+1}$.
    \item \textbf{Even Rows $2 \le r \le 2t-2$ and Rows $2t \le r \le n-1$:} $S_R(r, 2t+3) = 2t - 2n - 3 \equiv 2nk - 2n + 2t - 2 \pmod{2nk+1}$.
\end{itemize}

\begin{lem}\label{lem:rowaddedsums_Dh}
With the preceding partial sum $S_R(r, 2t+3)$ given above, for any row $0 \le r \le n-1$, summing the entries of diagonals $\mathbf{D_{n/2+1+2i}}$ and $\mathbf{D_{n/2+2i}}$ left-to-right starting in cell $(r, r+n/2-2t+1)$ with $0 \le i \le t-1$ yields the next $2t$ partial sums $S_R(r, 2t+4), \dots, S_R(r, 4t+3)$. The pair of partial sums modulo $2nk+1$ generated at step $i$ is determined by four cases based on the row $r$:

\begin{description}
    \item[Case 1: Row $0$.]
    The pair of sums is: $6n+4nt-6t+1-(4n-6)i \quad \text{and} \quad 2t-2i-5$ (which evaluates to $2nk$ for $i=t-2$ and $2nk-2$ for $i=t-1$).
    
    \item[Case 2: Odd Rows $1 \le r \le 2t-1$.]
    For $i < t-(r+1)/2$, the pair is $4nt-6t+4r+1-(4n-6)i$ and $2nk-6n+2t-2i-4$. For $i \ge t-(r+1)/2$, the pair is $4nt-6t+4r+1-(4n-6)i$ and $2nk-2n+2t-2i-4$.
    
    \item[Case 3: Even Rows $2 \le r \le 2t-2$.]
    For $i < t-r/2$, the pair is $4n+4nt-6t+4r+1-(4n-6)i$ and $2nk-2n+2t-2i-4$. For $i \ge t-r/2$, the pair is $4nt-6t+4r+1-(4n-6)i$ and $2nk-2n+2t-2i-4$.
    
    \item[Case 4: Rows $2t \le r \le n-1$.]
    The pair of sums is: $4nt-6t+4r+1-(4n-6)i \quad \text{and} \quad 2nk-2n+2t-2i-4$.
\end{description}
\end{lem}

Using these lemmas, we construct the complete partial sum sequence $S_R(r,p)$. We evaluate the first $k-2$ terms of $S_R(r,p)$ by sequentially adding: the $2t$ entries from diagonals $\mathbf{D_{2+2i}}$ and $\mathbf{D_{1+2i}}$ from \Cref{lem:rowaddedsums_D12}, the four consecutive base $H(n;6)$ array entries detailed in \Cref{tab:Hn6rowentries}, and the $2t$ entries from diagonals $\mathbf{D_{n/2+1+2i}}$ and $\mathbf{D_{n/2+2i}}$ from \Cref{lem:rowaddedsums_Dh}. The structure of the base $H(n;6)$ array and the starting cells of the diagonals require us to partition these sums into seven cases depending on row $r$: row $0$, odd and even rows for $1 \le r \le 2t-1$, odd rows for $2t+1 \le r \le n-5$, even rows for $2t \le r \le n-2$, and rows $n-3$ and $n-1$.

\begin{lem}\label{lem:rowsums_final}
Recall from \Cref{tab:Hn6rowentries} that the sum of the final two consecutive base array entries is $3$ for row 0 and $2n+3$ for all rows $1 \le r \le n-1$. By \Cref{lem:rowaddedsums_Dh}, the preceding partial sum $S_R(r, k-3)$ evaluates to $2nk-2$ for row 0, and $2nk-2n-2$ for all other rows. Adding the first of the two remaining base array entries yields the penultimate partial sum $S_R(r, k-2)$, which is determined by the row $r$:

\begin{description}
    \item[Row $0$:] $2nk-4n-1$
    \item[Rows $1 \le r \le n-1$:] $2nk-4n-2r-1$
\end{description}

In all cases, adding the final base array entry brings the cumulative sum to $S_R(r, k-1) = 0$, completing the row sum.
\end{lem}

\begin{thm}\label{thm:row_distinct}
For an $H(n;k)$ constructed as in \Cref{sec:construction}, the $k$ partial sums $S_R(r, p)$ in any row $r$ are distinct modulo $2nk+1$.
\end{thm}

\begin{proof}
By applying \Cref{lem:rowaddedsums_D12}, \Cref{lem:rowaddedsums_Dh}, and \Cref{lem:rowsums_final}, we can order the partial sums for each row into a strictly increasing sequence modulo $2nk+1$. Since the row affects the values of the added diagonals, we partition the rows into seven cases: row $0$, odd and even rows for $1 \le r \le 2t-1$, odd rows for $2t+1 \le r \le n-5$, even rows for $2t \le r \le n-2$, and rows $n-3$ and $n-1$.

\textbf{Row 0:}
\begin{align*}
  0 &< 1 < 2 < \dots < 2t-2 < 2t < n + 2t - 1 \\
  &< 2n + 2t - 1 < 10n - 5 < \dots < 4nt + 6n - 6t + 1 \\
  &< 2nk - 4nt - 6n + 8t - 2 < \dots < 2nk - 10n + 2t + 4 \\
  &< 2nk - 4n - 1 < 2nk - 2n + 2t - 1 < 2nk - 2 < 2nk
\end{align*}

\textbf{Odd Rows $1 \le r \le 2t-1$:}
\begin{align*}
  0 & < 2 < \dots < 2t - r - 1 < 4n + 4r - 5 < \dots < 4nt + 4r - 6t + 1 \\
  &< 2nk - 4nt - 6n - 4r + 8t - 2 < \dots < 2nk - 10n - 4r + 2t + 4 \\
  &< 2nk - 6n - 2 < \dots < 2nk - 6n + 2t - 4 \\
  &< 2nk - 6n + 2t - 2 < 2nk - 4n - 2r - 1 \\
  &< 2nk - 4n + 3 < \dots < 2nk - 4n + 2t + 1 \\
  &< 2nk - 3n + r + 2t + 3 < 2nk - 2n - 2 < \dots < 2nk - 2n + 2t - 4 < 2nk - n + 2t
\end{align*}

\textit{Note: The partial sum $2nk-6n-r+2t-2$ appears in this set of odd rows. Recall $t \ge 1$. While $2nk-6n-2 < 2nk-6n-r+2t-2$ is always true for $1 \le r \le 2t-1$, the exact position of $2nk-6n-r+2t-2$ relative to the sequence $2nk-6n-2 < \dots < 2nk-6n+2t-4$ changes depending on $r$. However, because $r$ is odd, $2nk-6n-r+2t-2$ is always odd. Since all terms in the sequence $2nk-6n+2i-2$ are even, this particular partial sum is guaranteed to be distinct no matter where it falls in this range.}

\textbf{Even Rows $2 \le r \le 2t-2$:}
\begin{align*}
  0 & < 2 < \dots < 2t < r + 2t - 1 \\
  &< 2r + 2t - 1 < 4n + 4r - 5 < \dots < 4nt + 4n + 4r - 6t + 1 \\
  &< 2nk - 4nt - 6n - 4r + 8t - 2 < \dots < 2nk - 6n - 4r + 2t + 4 \\
  &< 2nk - 4n - 2r - 1 < 2nk - 4n + r + 2t - 1 \\
  &< 2nk - 2n - 2 < \dots < 2nk - 2n + 2t - 4 < 2nk - 2n + 2t - 2
\end{align*}

\textbf{Odd Rows $2t+1 \le r \le n-5$:}
\begin{align*}
  0 & < 2 < \dots < 2t < n + r + 2t + 2 \\
  &< 3n + 2t - 1 < 4n + 4r - 5 < \dots < 4nt + 4r - 6t + 1 \\
  &< 2nk - 4nt - 2n - 4r + 8t - 2 < \dots < 2nk - 6n - 4r + 2t + 4 \\
  &< 2nk - 4n - 2r - 1 < 2nk - 2n - r + 2t - 2 \\
  &< 2nk - 2n - 2 < \dots < 2nk - 2n + 2t - 4 < 2nk - 2n + 2t - 2
\end{align*}

\textbf{Even Rows $2t \le r \le n-2$:}
\begin{align*}
  0 & < 2 < \dots < 2t < r + 2t - 1 \\
  &< 2r + 2t - 1 < 4n + 4r - 5 < \dots < 4nt + 4r - 6t + 1 \\
  &< 2nk - 4nt - 2n - 4r + 8t - 2 < \dots < 2nk - 6n - 4r + 2t + 4 \\
  &< 2nk - 4n - 2r - 1 < 2nk - 4n + r + 2t - 1 \\
  &< 2nk - 2n - 2 < \dots < 2nk - 2n + 2t - 4 < 2nk - 2n + 2t - 2
\end{align*}

\textbf{Row $n-3$:}
\begin{align*}
  0 & < 2 < \dots < 2t < n + 2t - 1 \\
  &< 3n + 2t - 1 < 8n - 17 < \dots < 4nt + 4n - 6t - 11 \\
  &< 2nk - 4nt - 6n + 8t + 10 < \dots < 2nk - 10n + 2t + 16 \\
  &< 2nk - 6n + 5 < 2nk - 3n + 2t + 1 < 2nk - 2n - 2 < \dots < 2nk - 2n + 2t - 4 \\
  &< 2nk - 2n + 2t - 2
\end{align*}

\textbf{Row $n-1$:}
\begin{align*}
  0 & < 2 < \dots < 2t < n + 2t + 1 \\
  &< 3n + 2t - 1 < 8n - 9 < \dots < 4nt + 4n - 6t - 3 \\
  &< 2nk - 4nt - 6n + 8t + 2 < \dots < 2nk - 10n + 2t + 8 \\
  &< 2nk - 6n + 1 < 2nk - 2n - 2 < \dots < 2nk - 2n + 2t - 4 < 2nk - 2n + 2t - 2 \\
  &< 2nk - 2n + 2t - 1
\end{align*}

Because the inequalities in the proof of \Cref{thm:row_distinct} arrange the $k$ partial sums of each row into a strictly increasing sequence within the interval $[0, 2nk]$, no two partial sums in the same row can be equivalent. Consequently, the partial sums for every row are pairwise distinct. Since we also established in \Cref{lem:rowsums_final} that the final sum of every row evaluates to $0 \pmod{2nk+1}$, these two facts together rigorously conclude the proof that the rows of the constructed array are simple. 
\end{proof}

Together with \Cref{thm:col_distinct}, this shows the Heffter arrays constructed in this manner are globally simple.

\section{Conclusion}

In this paper, we constructed globally simple integer Heffter arrays $H(n;k)$ for all even $n$ and $k \equiv 2 \pmod{4}$ ($k \ge 6$), thereby closing the last major open problem in the existence spectrum for globally simple square integer Heffter arrays. While the square case is now well understood, results on globally simple rectangular Heffter arrays $H(m,n;s,k)$ remain sparse \cite{dinitz2024}. Extending globally simple constructions to the rectangular case is an active area of research that may provide fruitful results.

\paragraph{Acknowledgement:} The authors would like to acknowledge support from TÜBİTAK 1001 grant number: 124F360.

%%%%%%%%%%%%%%%%%%%%%%%%%%%%%%%%%%%%%%%%%%%%%%%%%%%%%%%%%%%%%%%%%%%%%%%%%%%%%%%%%%%%%%%%%%%%%%%%%%%%%%%%%%%%%%%%%%%%%%
\bibliography{references_k_2_v5}

\begin{thebibliography}{15}
\providecommand{\natexlab}[1]{#1}
\providecommand{\url}[1]{\texttt{#1}}
\expandafter\ifx\csname urlstyle\endcsname\relax
  \providecommand{\doi}[1]{doi: #1}\else
  \providecommand{\doi}{doi: \begingroup \urlstyle{rm}\Url}\fi

\bibitem[Archdeacon(2015)]{archdeacon2015}
D.S. Archdeacon.
\newblock Heffter arrays and biembedding graphs on surfaces.
\newblock \emph{Electron. J. Combin.}, 22:\penalty0 5--35, 2015.

\bibitem[Archdeacon et~al.(2015)Archdeacon, Dinitz, Donovan, and
  Yazıcı]{archdeacon2_2015}
D.S. Archdeacon, J.H. Dinitz, D.M. Donovan, and E.Ş. Yazıcı.
\newblock Square integer {Heffter} arrays with empty cells.
\newblock \emph{Des. Codes Cryptogr.}, 77:\penalty0 409--426, 2015.

\bibitem[Archdeacon et~al.(2016)Archdeacon, Dinitz, Mattern, and
  Stinson]{archdeacon2016}
D.S. Archdeacon, J.H. Dinitz, A.~Mattern, and D.R. Stinson.
\newblock On partial sums in cyclic groups.
\newblock \emph{JCMCC}, 98:\penalty0 327--342, 2016.

\bibitem[Bode and Harborth(2005)]{bode2005}
J.P. Bode and H.~Harborth.
\newblock Directed paths of diagonals within polygons.
\newblock \emph{Des. Math.}, 299:\penalty0 3--10, 2005.

\bibitem[Burrage et~al.(2020)Burrage, Cavenagh, Donovan, and
  Yazıcı]{burrage2020}
K.~Burrage, N.J. Cavenagh, D.M. Donovan, and E.Ş. Yazıcı.
\newblock Globally simple {Heffter} arrays {$H(n;k)$} when $k \equiv 0,3
  \pmod{4}$.
\newblock \emph{Discrete Math.}, 343:\penalty0 111787, 2020.

\bibitem[Caro and Yuster(1999)]{caro1999}
Y.~Caro and R.~Yuster.
\newblock Orthogonal decompositions and packing of complete graphs.
\newblock \emph{J. Combin. Theory Ser. A}, 88:\penalty0 93--111, 1999.

\bibitem[Cavenagh et~al.(2019)Cavenagh, Dinitz, Donovan, and
  Yazıcı]{cavenagh2019}
N.J. Cavenagh, J.H. Dinitz, D.M. Donovan, and E.Ş. Yazıcı.
\newblock The existence of square non-integer {Heffter} arrays.
\newblock \emph{Ars Math. Contemp.}, 17:\penalty0 369--395, 2019.

\bibitem[Cavenagh et~al.(2020)Cavenagh, Donovan, and Yazıcı]{cavenagh2020}
N.J. Cavenagh, D.M. Donovan, and E.Ş. Yazıcı.
\newblock Biembeddings of cycle systems using integer {Heffter} arrays.
\newblock \emph{J. Combin. Des.}, 28:\penalty0 900--922, 2020.

\bibitem[Costa et~al.(2018)Costa, Morini, Pasotti, and Pellegrini]{costa2018}
S.~Costa, F.~Morini, A.~Pasotti, and M.A. Pellegrini.
\newblock Globally simple {Heffter} arrays and orthogonal cyclic cycle
  decompositions.
\newblock \emph{Austral. J. Combin.}, 72:\penalty0 549--593, 2018.

\bibitem[Costa et~al.(2026)Costa, Della~Fiore, Fontana, and Vena]{costa2026}
S.~Costa, S.~Della~Fiore, M.~Fontana, and L.~Vena.
\newblock Graham conjecture on small sets in abelian groups, 2026.
\newblock URL \url{https://arxiv.org/abs/2603.20961}.

\bibitem[Dinitz and Pasotti(2024)]{dinitz2024}
J.H. Dinitz and A.~Pasotti.
\newblock A survey of {Heffter} arrays.
\newblock In C.J. Colbourn, editor, \emph{New Advances in Designs, Codes and
  Cryptography}, volume~86, pages 353--392. Springer Nature Switzerland, 2024.

\bibitem[Dinitz and Wanless(2017)]{dinitz2017}
J.H. Dinitz and I.M. Wanless.
\newblock The existence of square integer {Heffter} arrays.
\newblock \emph{Ars. Math. Contemp.}, 13:\penalty0 81--93, 2017.

\bibitem[Gordon(1961)]{gordon1961}
B.~Gordon.
\newblock Sequences in groups with distinct partial products.
\newblock \emph{Pacific J. Math.}, 11:\penalty0 1309--1313, 1961.

\bibitem[Heffter(1896)]{heffter1896}
L.~Heffter.
\newblock Über nachbarconfigurationen, triplesysteme und metacyklische
  gruppen.
\newblock \emph{Deutsche Mathem. Vereinig. Jahresber.}, 5:\penalty0 67--69,
  1896.

\bibitem[Pelttari et~al.(2026)Pelttari, Küçükçifçi, and
  Yazıcı]{pelttari2026}
E.~Pelttari, S.~Küçükçifçi, and E.Ş. Yazıcı.
\newblock Globally simple {Heffter} arrays {$H(n;k)$} with $k \equiv 1
  \pmod{4}$, 2026.
\newblock URL \url{https://arxiv.org/abs/2604.20252}.

\end{thebibliography}

\end{document}